\documentclass[10pt,journal]{IEEEtran}

\ifCLASSINFOpdf
   \usepackage[pdftex]{graphicx}
    \usepackage{epstopdf}
\else
   \usepackage[dvips]{graphicx}
\fi

\usepackage{amssymb}
\usepackage{amsmath}
\usepackage{alltt}
\usepackage{framed}
\usepackage{epstopdf}

\begin{document}

\title{Enhancing the Performance and Robustness of the FEAST Eigensolver}

\author{\IEEEauthorblockN{Brendan Gavin and Eric Polizzi}\\
\IEEEauthorblockA{Department of Electrical and Computer Engineering, \\
University of Massachusetts - Amherst, MA 01003, USA \\
E-mail: bgavin@umass.edu, polizzi@ecs.umass.edu}
}


\maketitle

\begin{abstract}
The FEAST algorithm is a subspace iteration method that uses a spectral projector as a rational filter 
in order to efficiently solve interior eigenvalue problems in parallel.
Although the solutions from the FEAST algorithm converge rapidly 
in many cases, convergence can be slow in situations where the eigenvalues of a matrix are densely populated 
near the edges of the search interval of interest, which can be detrimental
 to parallel load balancing. This work introduces two methods that allow one to improve the convergence robustness 
of the FEAST algorithm in these situations without having to increase the amount of computation.
Selected numerical examples are presented and discussed.

\end{abstract}
\IEEEpeerreviewmaketitle
\section{Introduction}

FEAST \cite{polizziFEAST,tangFEAST} is a subspace iteration algorithm
 for solving eigenvalue problems 

\begin{equation}
Ax=\lambda Bx,\ \  A \in \mathbb{R}^{n \times n},
\end{equation}

\noindent by finding the eigenvectors $x$ whose
eigenvalues $\lambda$ lie in some interval ${\cal I}=(\lambda _{min},\lambda _{max} )$ of the user's choosing. 
In this paper we consider the Hermitian standard eigenvalue problem for simplicity (i.e. $A=A^H$ and $B\equiv I$), 
but the FEAST algorithm can be extended straight-forwardly to the generalized and non-Hermitian eigenvalue problems
 as well \cite{kpt16}.

 FEAST belongs to the broader family of 
contour integration eigensolvers \cite{sakurai1,sakurai2,sakurai3,austin1}, but it can also be accurately described as an optimal subspace iteration procedure.
A conventional subspace iteration consists of multiplying a trial subspace by the matrix
 $A$ and then orthogonalizing it with the Rayleigh-Ritz procedure; this process is repeated 
iteratively until the subspace converges. The FEAST algorithm operates similarly, 
but rather than multiplying the trial subspace by $A$, one instead multiplies the trial subspace by the spectral 
projector matrix $\rho (A)$. The matrix $\rho (A)$ is given by the complex contour integral
\begin{align}
\rho (A) &= \frac{1}{2\pi i}\oint_{\cal C} (zI-A)^{-1}dz,\  \label{contour_int}
\end{align}
\noindent where $\cal C$ is a closed contour in the complex plain that exactly encloses the 
interval $\cal I$ on the real axis. 
The function $\rho (\lambda)$ applied to a real number 
$\lambda$ is a filter function that returns 1 when $\lambda \in \cal I$ and 0 otherwise. 
As a result, the matrix $\rho (A)$ is a spectral projector whose image is the subspace
 that is spanned by only the eigenvectors of $A$ whose eigenvalues lie in $\cal I$.
 Multiplication of a vector by $\rho(A)$ projects it into that subspace, and in this
 way the FEAST algorithm finds only the eigenvectors whose eigenvalues lie in
 $\cal I$ \cite{tangFEAST}. The algorithm is outlined in Appendix.
The contour integral (\ref{contour_int}) has no
 general analytical expression, so in practice the multiplication of a matrix 
$X$ by $\rho (A)$ is approximated by using some numerical integration quadrature rule
\begin{align}
\rho (A)X = \frac{1}{2\pi i}\oint_{\cal C} (zI-A)^{-1}Xdz\ \approx \sum_{i=1}^{n_c} \omega _i (z_iI-A)^{-1}X, \label{contour_int_quad}
\end{align}
\noindent where $n_c$ is the number of quadrature points, and each term $(z_iI-A)^{-1}X$ is found by using a 
linear system solver with
the column vectors of the matrix $X$ as the right hand sides of the linear system. 
A variety of quadrature rules are possible; in this work we use Gauss quadrature. 

The benefits of using FEAST over a traditional subspace iteration technique are twofold. The first benefit 
is that, by finding only the eigenvectors whose eigenvalues lie in a certain interval, eigenvalue 
problems can be solved in parallel by solving for the eigenvector/eigenvalue pairs in different intervals independently. 

The second benefit has to do with the rate of convergence. In a subspace iteration algorithm operating on 
a subspace of dimension $m_0$, the eigenvector with the $i^\text{th}$ largest eigenvalue magnitude $|\lambda_i|$ 
 converges at a rate of $|\lambda_i|/|\lambda_{m_0+1}|$, where $\lambda_{m_0+1}$ is the eigenvalue with the
 $(m_0+1)^\text{th}$ largest magnitude. For a typical subspace iteration this means that the rate of 
convergence depends strongly on the eigenspectrum of $A$. By using FEAST, the rate of convergence 
becomes $\rho(\lambda_i)/\rho(\lambda_{m_0+1})$ \cite{tangFEAST}, where $\lambda_i$ is now the eigenvalue
 with the $i^\text{th}$ largest value of $\rho(\lambda)$, and $\lambda_{m_0+1}$ is the eigenvalue with 
the $(m_0+1)^\text{th}$ largest value of $\rho(\lambda)$. 

The ratio $\rho(\lambda_i)/\rho(\lambda_{m_0+1})$ can be made arbitrarily large by either increasing the 
accuracy of the quadrature rule 
(\ref{contour_int_quad}) by increasing $n_c$, or by increasing the size of the subspace $m_0$; the eigenpairs
anywhere in the spectrum can thus be found rapidly. It is not uncommon to be able to achieve
a convergence rate of $10^4$ with $n_c=8$ and a subspace size of $m_0\approx1.5m$, where $m$ is the exact
 number of eigenvalues that lie in the interval~$\cal I$.
Because of these remarkable convergence properties, as well as its robustness and its ability to exploit parallelism at multiple levels, the FEAST algorithm and associated software package ({\em www.feast-solver.org}) have been very well received by the HPC community. The FEAST algorithm is currently featured as the principle HPC eigenvalue solver in the Intel Math Kernel Library (MKL).

\begin{figure}[h!]
\centering
\includegraphics[width=\linewidth]{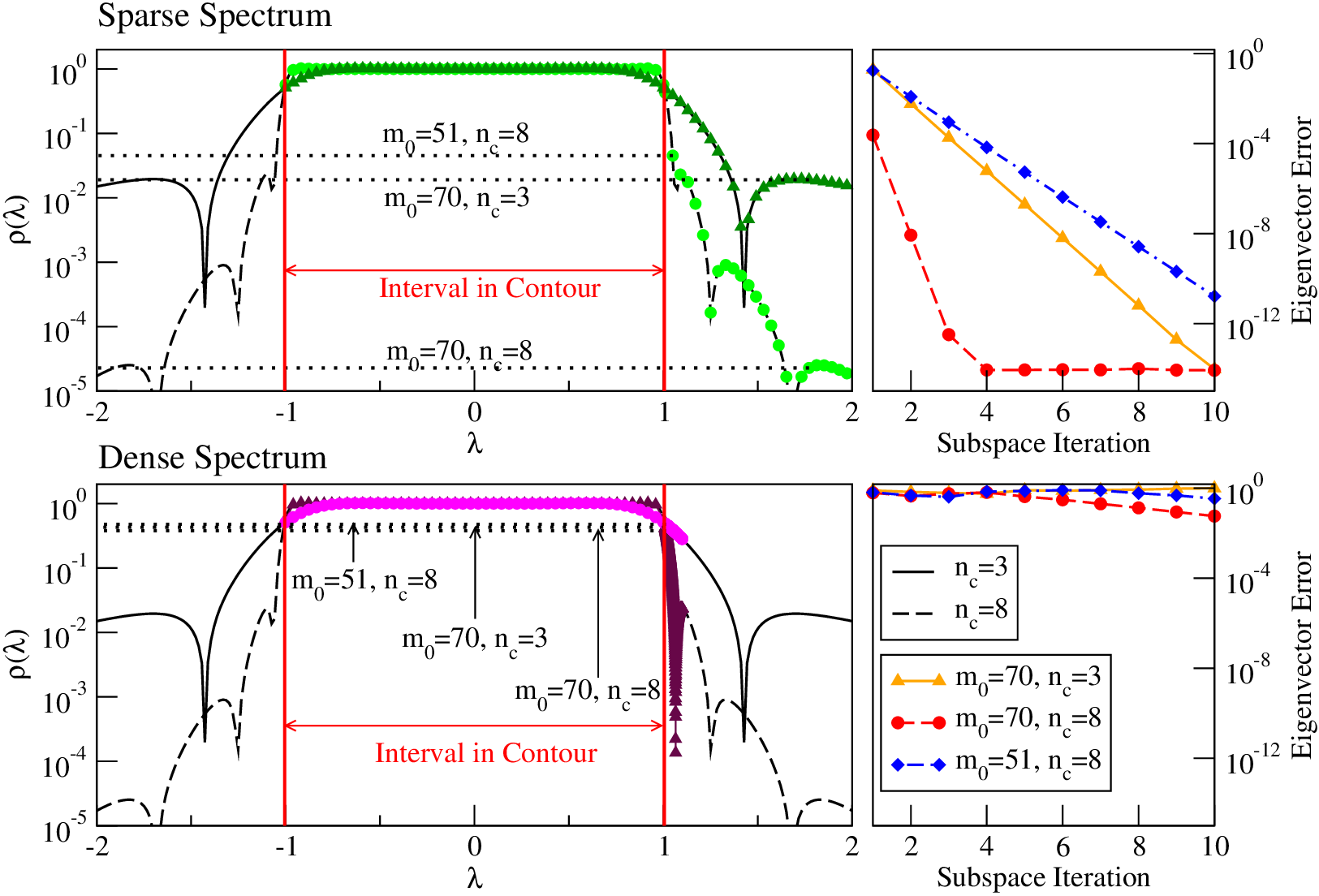}
\caption{Two test cases illustrating the difference in the convergence rate of FEAST for a matrix with a sparsely 
packed eigenvalue spectrum outside the contour interval (top plots) and a matrix with a densely packed eigenvalue
 spectrum outside the contour interval (bottom plots). Each matrix is dimension 545 with 50 eigenvalues inside of
 the contour interval and 495 eigenvalues outside of the contour interval. The plots on the right show the convergence
 of the maximum eigenvector error for various values of the parameters $m_0$ and $n_c$. The plots on the left show 
the value of $\rho(\lambda)$ for $n_c=$ 3 and $n_c= $8 plotted with solid and dashed curves, with the locations of 
the eigenvalues of the matrix indicated by plot markers. The locations of $\lambda_{m0+1}$ and the values of 
$\rho(\lambda_{m0+1})$ are indicated with dotted horizontal lines for the same several values of $m_0$ and $n_c$. 
The matrix with the sparsely packed spectrum converges well, whereas the matrix with densely packed spectrum barely
 converges at all.
}
\label{rfplot}
\end{figure}

The convergence rate of FEAST is not entirely insensitive to the spectrum of $A$, however. 
In situations where the eigenvalues of $A$ are packed many times more closely together immediately outside 
of $\cal I$ than they are inside of $\cal I$, the rate of convergence can be very slow. This is illustrated
 in Figure \ref{rfplot}. The top plots in Figure \ref{rfplot} illustrate the situation where the density 
of the eigenvalue spectrum is the same both inside and outside the interval $\cal I$, and the bottom plots
 in Figure \ref{rfplot} illustrate the situation where the density is much larger outside of the interval 
$\cal I$ than it is inside of the interval $\cal I$. The error at each FEAST subspace iteration is plotted 
for several values of $n_c$ and $m_0$, and the corresponding values of $\lambda_{m_0+1}$ and $\rho(\lambda_{m_0+1})$ 
are indicated with horizontal dotted lines in the plots on the left in order to illustrate the effects of these 
parameters on convergence for both the dense spectrum and the sparse spectrum.

Even in situations that are less pathological than the one illustrated in the bottom plot of Figure \ref{rfplot}, 
the varying density of the spectrum of $A$ can have negative implications for parallel load balancing.
 We can find the eigenpairs of $A$ in parallel by dividing the spectrum of $A$ into several non-intersecting
 intervals and then solving the eigenvalue problem for each interval separately and in parallel. When we do 
this, some intervals may converge more quickly than others due to the varying density of the spectrum, even
 if every interval contains the same number of eigenvalues. It is possible to speed up convergence in a given 
interval by increasing $n_c$ or $m_0$, but this does not reduce the amount of computation required; increasing
 $n_c$ or $m_0$ increases the number of linear systems that need to be solved with each iteration, and the
 solution of the linear systems for the quadrature rule in equation (\ref{contour_int_quad}) is where most
 of the computation in the FEAST algorithm occurs. 

We would ideally like to be able to use parallel resources as efficiently as possible, performing the same amount 
of computation for each interval in which we solve the eigenvalue problem. We therefore would like to
improve the convergence rate of FEAST in situations where the spectrum of $A$ results in slow or varying convergence
 rates, but without having to solve additional linear systems in order to do so.

In Ref. \cite{zolotarev}, this problem is addressed with the introduction of the Zolotarev quadrature 
that produces a very steep slope for the filter at the interval endpoints, which then leads to the same convergence 
rate between different contours. The Zolotarev approach presents, however, few limitations. The first limitation is that the convergence rate is fixed and cannot be improved while increasing $m_0$, and it will thus underperform in comparison
with Gauss quadrature, for example, in situations where the spectrum is sparsely packed or uniformely distributed 
 (e.g. top plot of Figure \ref{rfplot}); 
The second limitation is that the Zolotarev approach cannot be extended to the non-Hermitian problem where the eigenvalues are located in the complex plane.

In this work we propose a more general set of alternatives that use ``accelerated subspace approach'' strategies in order to improve the convergence robustness of FEAST regardless of which quadrature rule is being used.

\section{Accelerating the FEAST Subspace Iteration}

Previous research \cite{gavinFEAST} and the observation that larger subspace sizes $m_0$ increase the 
rate of convergence for FEAST suggest that we may be able to improve convergence by finding ways to increase
 the size of the subspace that is used in the Rayleigh-Ritz procedure. If we can do this without having
 to solve additional linear system right hand sides when performing the numerical quadrature in equation
 (\ref{contour_int_quad}), then we may improve the convergence rate of FEAST without having to do too 
much additional computation. In the following subsections we discuss two ways of expanding the FEAST
 subspace size without solving additional linear systems.

\subsection{Method 1: Expand Subspace Using Previous Subspaces}
\label{xfeast_sec}
In a typical FEAST subspace iteration, the trial subspace $X_i$ from the previous iteration 
is discarded and replaced with the filtered  subspace $\rho (A)X_i$ (Step 1, FEAST Algorithm). 
Rather than discarding the previous subspace $X_i$, we might instead append the new, filtered 
subspace to the old one before performing the Rayleigh-Ritz procedure; by doing this we can 
increase the dimension of the subspace by $m_0$ without having to solve additional linear 
systems. Step 1 of FEAST might then look like this:

\medskip
\begin{minipage}{0.8\linewidth}
{\footnotesize
\begin{description}
\item[1.] Filter the trial subspace and append it to the columns of the old one: 
\begin{center} $X'=[X_i\ \ \rho (A)X_i]$\end{center}
\end{description}
}
\end{minipage}
\medskip

\noindent where we form $X'$ by appending the column vectors of $\rho (A)X_i$ to the matrix for the previous
 subspace $X_i$. We could repeat this process several times in order to build up a total subspace size of
 $s \times m_0$, after which we could keep the subspace size constant by removing old subspaces before 
adding new ones at each subspace iteration. 

It would not be surprising if this modification of FEAST were to improve its the convergence rate; by
 expanding the subspace in this way, we are essentially building a Krylov subspace technique wherein
 we multiply our prospective subspace by powers of $\rho(A)$ rather than by powers of $A$.  
 
If we modify step 1 of FEAST in this way then we have to make a few other modifications as well. 
Step 2(i) of FEAST requires the solution of the reduced eigenvalue problem $A' q=\lambda B' q$, 
so we need to ensure that $B'=X'^TX'$ is symmetric positive definite; because we append the 
filtered subspace to the old one, this is no longer guaranteed. We therefore need to add 
another step to FEAST: orthogonalize the matrix $X'$ before doing the Rayleigh-Ritz procedure. 
This can be done by using the QR decomposition or the singular value decomposition (SVD) of $X'$. 
For the research presented here, we orthogonalize $X'$ by taking its SVD and setting $X'$ equal
 to the left singular vectors: 
\begin{equation}
X'=U\Sigma V^T\ \rightarrow \ X'=U.
\end{equation}
\noindent In particular, we do this by diagonalizing $X'^TX'$
\begin{equation}
X'^TX'=V\Sigma ^2 V^T \rightarrow U=X'V\Sigma^{-1},
\end{equation}
\noindent and then retaining the first $m_0$ columns of $U$. 
Although this is less numerically stable than QR, we have found that it offers performance benefits in terms of the speed of the 
orthogonalization.

We call this algorithm 'expanding subspace FEAST'; see the XFEAST algorithm in Appendix.

The implementation of XFEAST we use in this paper also involves expanding the subspace to its full
 size before doing the first Rayleigh-Ritz procedure. The subspace size can be increased incrementally,
 with the Rayleigh-Ritz procedure being done in between each subspace expansion, but there is no reason
 to do this unless one expects that the algorithm might converge before the subspace size has reached its limit. 

%
%
%
%

The part of step {\bf 2(ii)} of XFEAST that specifies that one must select the desired eigenvectors is
 required because the subspace $X$ is expanded beyond just the filtered subspace. In conventional FEAST
 iterations the Rayleigh-Ritz procedure will find all of the $m$ eigenvectors whose eigenvalues are in
 the interval $\cal I$, plus the $m_0-m$ eigenvectors whose eigenvalues  are closest to, but still
 outside of, $\cal I$. In XFEAST, because the subspace is expanded beyond the size $m_0$, the
 Rayleigh-Ritz procedure will find all of those $m_0$ eigenpairs plus many more. Due to
 numerical errors it may even find eigenpairs for the Rayleigh-Ritz matrix $A'$ that do 
not exist for the original matrix $A$.

Since step {\bf 4(i)} of XFEAST requires a subspace of dimension $m_0$ to filter with $\rho(A)$ 
for the next iteration, we must select $m_0$ of the $s\times m_0$ eigenpairs that are produced 
by step {\bf 2(i)}. Here, we are using a few steps of sorting. 
First, we calculate the error residuals for all $s\times m_0$ eigenpairs from {\bf 2(i)}. 
We then select all of the eigenpairs whose eigenvalues lie inside the interval 
${\cal I}=(\lambda _{min}, \lambda _{max})$.
If fewer than $m_0$ eigenpairs are found whose eigenvalues lie inside the contour interval,
 we then select additional eigenpairs from outside the contour as well, preferentially
 selecting those eigenpairs with the lowest residuals. 


\subsection{Method 2: Expand Subspace Using Eigenvector Residuals}

The other piece of information that the typical FEAST iteration generates (and which is otherwise discarded) 
is the eigenvector error residuals. Step 3 of the normal FEAST algorithm computes the eigenvector
 error residuals $r_k=Ax_k-\lambda _k x_k$, and uses the one with the largest norm as a measure of
 the accuracy of the current subspace estimate. 

Because the current eigenpair estimates provided at each iteration of FEAST come from the 
Rayleigh-Ritz procedure, 
the inner product of any of the estimated eigenvectors with any of the residual vectors is zero:
$x_j^Tr_k=0,\ \forall \ 1\leq j,k\leq m_0$. One can show this by using the fact that $x_k=X'q_k$:
\begin{align}
\begin{split}
x_j^Tr_k &=x_j^TAx_k-x_j^T\lambda_k x_k \\
&=q_j^TX'^TAX'q_k-\lambda_k q_j^TX'^TX'q_k\\
&=\delta_{jk}\lambda_k - \lambda_k \delta_{jk}=0.
\end{split}
\end{align}
\noindent If $R$ is the matrix of column vectors $r_k$, then its column vectors span a subspace that 
is orthogonal to the current estimated solution subspace $X$. We can therefore perform another 
Rayleigh-Ritz procedure in the subspace spanned by the combined columns of $X$ and $R$ without 
having to orthogonalize the column vectors of $R$ with respect to those of $X$ in order to 
ensure that $X'^TX'$ is symmetric positive definite. This allows us to improve the estimated
 subspace without having to solve any additional linear systems and without having to do any
 orthogonalization procedure. A modified FEAST algorithm using this approach is given in the RFEAST algorithm in Appendix.

Again, it would not be surprising if adding eigenvector residuals to the subspace were to help improve the convergence rate; the eigenvalue algorithm LOBPCG \cite{knyazev2001toward} also works by including an eigenvector residual block in the search subspace. 

Step {\bf 2(ii)} of RFEAST again requires that we select the desired eigenpairs from amongst 
the ones produced by the Rayleigh-Ritz procedure. This is done in the same way as for XFEAST.

Measuring the error on the estimated subspace for RFEAST generally requires more care than in
 XFEAST or FEAST. The Rayleigh-Ritz procedure for RFEAST tends to produce many eigenpairs 
that do not exist in the spectrum of the full size matrix $A$ due to numerical error, and 
many of these spurious eigenpairs have eigenvalues that fall inside the interval $\cal I$. 

In order to return the correct estimated eigenpairs and estimate the error on them, we must
 select only the eigenpairs inside $\cal I$ that are not spurious. We do this by determining
 how many eigenpairs we should expect to find in that interval, and then taking that number
 of eigenpairs inside $\cal I$ with the lowest residuals to be the eigenpairs of interest. 
We determine the number of eigenpairs to expect by counting the number of eigenpairs found
 during the first Rayleigh-Ritz procedure of each subspace iteration (i.e. during iteration 
$j=1$ in step {\bf 2} of RFEAST); the subspace used for the first Rayleigh-Ritz procedure
 is just the conventional FEAST subspace, and so we will not yet have produced the 
proliferation of spurious eigenpairs that comes from expanding the subspace by 
using the eigenvector residuals.

\section{Results and Comparisons}
\label{results_section}

We demonstrate the convergence properties of these modified FEAST algorithms with several example matrices.

\begin{figure}[h!]
\begin{center}
\textbf{FEAST Iterations for Sparse and Dense Eigenspectra} 
\end{center}
\par
\medskip
\includegraphics[width=0.9\linewidth]{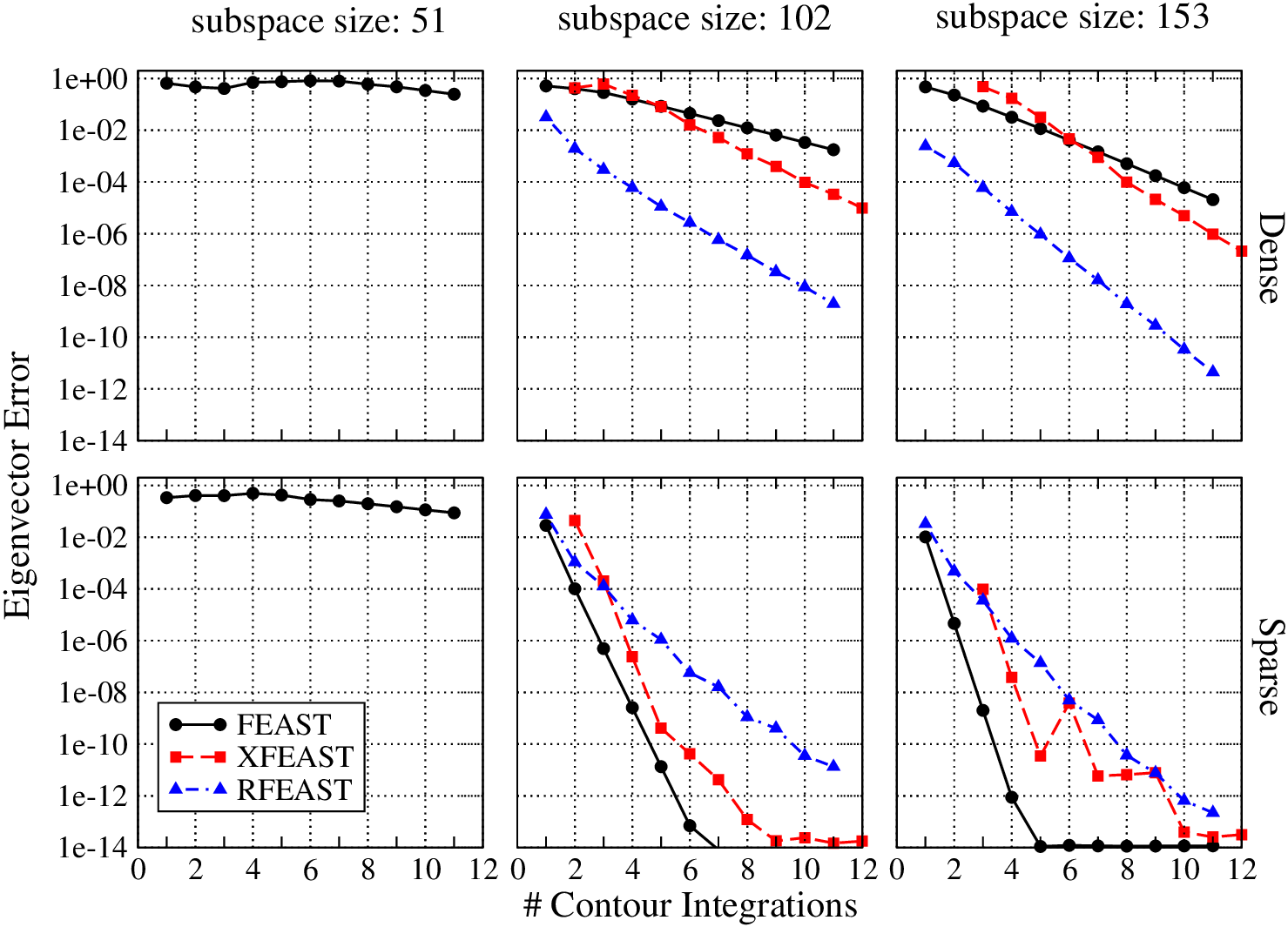}
\caption{Plots showing the eigenvector error residual versus number of contour integrations for each of 
the three FEAST variations for various subspace dimensions. The top row of plots show the results when 
using a matrix with a densely packed eigenspectrum outside of the FEAST interval, and the bottom row of
 plots show the results when using a matrix with a sparsely packed eigenspectrum outside of the FEAST 
interval. Both matrices are dimension 545, and we search for 50 eigenvalues. The number above each 
plot indicates the size of the subspace being used. 
The ``Dense'' results are for $n_c=8$ and the ``Sparse'' results are for $n_c=3$; convergence is too fast
 for good illustration with $n_c=8$ for the ``Sparse'' matrix. }
\label{rxf_vs_f_both}
\end{figure}

Figure \ref{rxf_vs_f_both} shows the eigenvector error residual at each subspace iteration of the FEAST,
 XFEAST, and RFEAST algorithms as applied to two different real symmetric matrices, for several different
 subspace sizes. Both matrices are dimension 545 and have the same eigenvectors, with 50 eigenvalues
 inside the FEAST interval $\cal I=$[-1,1]. 

One matrix, labeled ``Sparse'' in Figure \ref{rxf_vs_f_both}, has the other 495 eigenvalues in the interval
 [1.01, 20.81], whereas the one labeled ``Dense'' has those 495 eigenvalues in the interval [1.01, 1.1]. 
That is, the ``Sparse'' matrix has sparsely-packed eigenvalues outside of $\cal I$, and the ``Dense'' 
matrix has densely-packed eigenvalues outside of $\cal I$. 

For the FEAST iterations the value of $m_0$ is the same as the subspace size, whereas for the XFEAST
 and RFEAST iterations $m_0$ is always set at 51, and the full subspace size is generated by one or 
the other subspace expansion method. XFEAST and RFEAST thus solve the same number of linear systems
 for each contour integration, regardless of the subspace size, whereas FEAST solves more linear 
systems for larger subspace sizes.

Despite solving many fewer linear system right hand sides per iteration (the number of linear system 
right hand sides per iteration is $n_c\times m_0$), XFEAST and RFEAST outperform FEAST for the ``Dense''
 matrix on a per-subspace iteration basis. This is not the case for the ``Sparse'' matrix. Nonetheless,
 even for the ``Sparse'' matrix, XFEAST and RFEAST do a similar amount of total computation for a given 
level of convergence.

\begin{figure}[h!]
\begin{center}
\textbf{Eigenvector Residual vs. \# Linear RHS Solved} 
\end{center}
    \includegraphics[width=\linewidth]{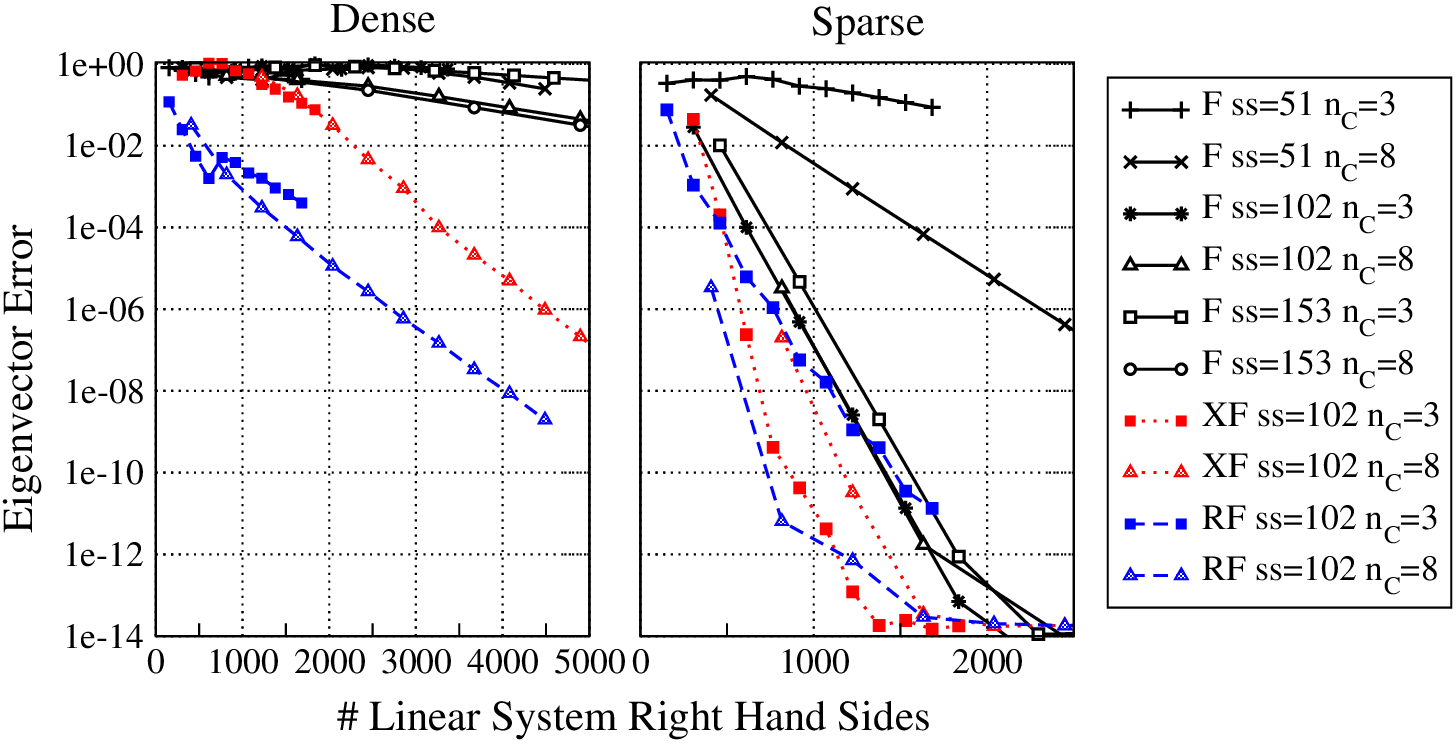}
\caption{Plots showing eigenvector residual versus the number of linear system right hand sides solved to
 reach that level of convergence, for 
both a matrix with a dense eigenspectrum outside the interval of interest and a matrix with a sparse 
eigenspectrum outside the interval of interest, 
using FEAST, XFEAST, and RFEAST for various subspace sizes $m_0$ and numbers of quadrature points $n_c$.
 XFEAST and RFEAST consistently require fewer linear system solutions than regular FEAST does in order 
to reach the same level of accuracy, with the difference being fairly dramatic in the dense eigenspectrum case.}
\label{rxf_vs_f_linsys}
\end{figure}

Figure \ref{rxf_vs_f_linsys} shows the amount of eigenvector error per number of linear system right 
hand sides solved for the same two matrices, for various values of $m_0$ and $n_c$. The advantages of
 XFEAST and RFEAST are especially clear here; the rate of convergence per linear system right hand side
 solved, which is the majority of the computation in the FEAST algorithm, depends primarily on which
 algorithm is used in the case of the ``Dense'' matrix, with XFEAST and RFEAST clearly outperforming FEAST.
 XFEAST and RFEAST also outperform FEAST for the ``Sparse'' spectrum matrix, but here the difference is less dramatic.

\begin{figure}[h!]
  \centering
\textbf{FEAST Applied to Electronic Structure Spectrum} 

\medskip

    \includegraphics[width=\linewidth]{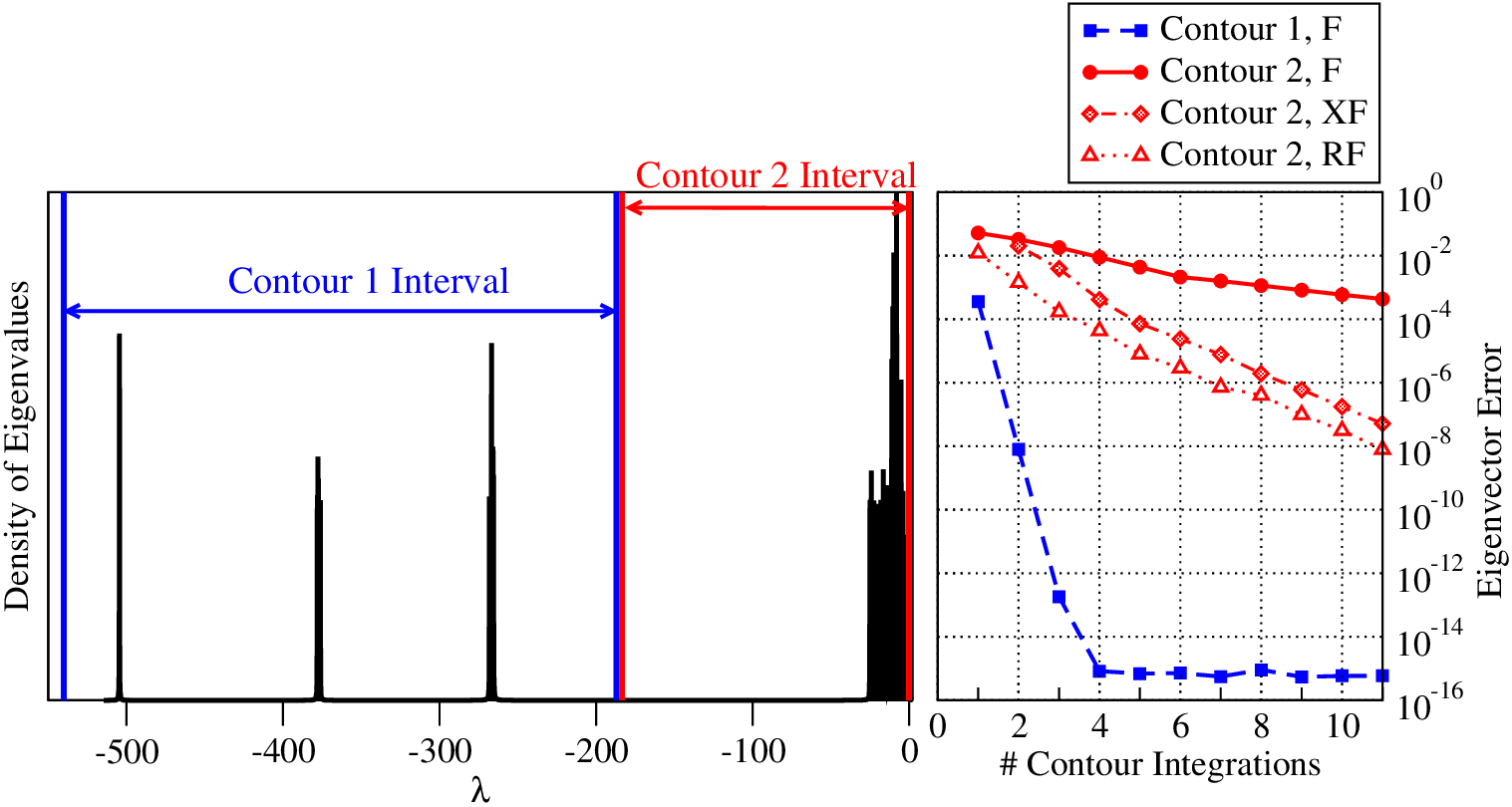}
\caption{Plots showing the application of the FEAST variations to a matrix derived from electronic structure
 theory. Left plot shows the density of the eigenspectrum of the matrix, divided into two intervals, and the right plot shows the 
convergence of the eigenvector
error for the various FEAST algorithms applied to the two intervals. The 14 eigenpairs in the "Contour 1 Interval" were calculated using a base subspace size of $m_0=17$, and the 43 eigenpairs in the "Contour 2 Interval" were calculated by using a base subspace size of $m_0=46$. Both the XFEAST and the RFEAST runs for "Contour Interval 2" use a total subspace size of $3m_0$, with the subspace having expanded twice by using either the previous FEAST iteration solutions or the eigenvector residuals.
}
\label{caf_cont}
\end{figure}

Figure \ref{caf_cont} illustrates the results of solving an eigenvalue problem whose spectrum derives 
from electronic structure calculations \cite{zolotarev,levin}. Our work here has been motivated by applications of this kind.

The left plot of Figure \ref{caf_cont} shows the density of the eigenspectrum of the Hamiltonian matrix
 for the ground state of of a Caffeine molecule. The density of the eigenspectrum shows several distinct
 peaks, and one potential partitioning of the spectrum into two intervals is shown with red and blue lines.
 The right plot of Figure \ref{caf_cont} shows the convergence trajectory when FEAST is used on each of these 
intervals separately. For the interval encompassed by Contour 2 (shown in red), the rightmost edge of which 
passes through a very dense region in the eigenspectrum, we also show the result of using XFEAST and RFEAST 
in order to try to achieve better convergence than is possible with FEAST.

The first interval (shown in blue), which has no eigenvalues immediately near its edges, converges rapidly, 
much like the second and third columns of the sparse example in Figure \ref{rxf_vs_f_both}.  The second interval, 
which has its upper limit passing through the middle of a dense group of eigenvalues, converges very slowly when
 using FEAST. This is the sort of problem that we seek to address. 

Using XFEAST and RFEAST, we can improve the final eigenvector error residual for the second, more challenging 
interval by more than four orders of magnitude. Still, this does not achieve ideal load balancing because the
 first interval has both of its edges in regions that are completely empty of eigenvalues, and so it converges 
very quickly. Better load balancing can only be achieved by dividing the spectrum in a less arbitrary way, which 
will require that we estimate the spectrum of a matrix before diagonalizing it. This is a subject of continuing 
research.

%

\section{Conclusion}

The results in Section \ref{results_section} show that we can indeed improve the convergence rate of FEAST 
without solving additional linear systems by expanding the FEAST subspace through other means. This is 
particularly helpful in situations where the spectrum of the matrix at hand makes convergence difficult. 
Doing so comes at the price of having to use additional memory to store the expanded subspace; when using
 enough parallelism (and therefore a large enough number of intervals), however, we expect that memory
 will not be a constraint because the initial size of the subspace for each interval can be made almost 
arbitrarily small. As the results in Figure \ref{caf_cont} show, though, this alone is not yet a fully 
satisfactory solution for achieving load balancing. 
Future work will consist of using this research to build on the efforts of others in order to estimate the 
eigenvalue distribution of a matrix \cite{stochasticEst} and efficiently divide the eigenvalue interval
 of interest \cite{kajpust}. We expect that, by combining our work here with these techniques for measuring
 and dividing and the eigenvalue spectrum of a matrix, we can achieve ideal load balancing in an automated 
way for arbitrary matrices.

\section*{Acknowledgments}
The authors wish to acknowledge helpful discussions
with Dr. Ping Tak Peter Tang and Dr. Yousef Saad.
This material is supported by NSF under Grant \#CCF-1510010.

%
%

\bibliographystyle{IEEEtran}
\bibliography{kfpap2col}

\appendix

\bigskip

\noindent \begin{minipage}{\linewidth}
\begin{framed}
{\footnotesize
\noindent \textbf{FEAST Algorithm} 

\noindent \hrulefill


\noindent {\bf Start with:} Matrix $A \in \mathbb{R}^{n \times n}$ to be diagonalized, 
interval ${\cal I}=(\lambda _{min},\lambda _{max} )$ wherein 
fewer than $m_0$ eigenvalues are expected to be found, 
initial guess $X_0 \in \mathbb{R}^{n\times m_0}$ for the subspace spanned by the solution to the eigenvalue problem
$Ax=\lambda x,\ \lambda \in {\cal I}$.
\begin{description}
\setlength{\itemsep}{3pt}
\item[1.] Filter the subspace $X_i$ to remove eigenvectors whose eigenvalues do not lie in the interval ${\cal I}$: 
$X'=\rho (A)X_i$
\item[2.] Perform Rayleigh-Ritz procedure to find a new estimate for eigenvalues and eigenvectors: 
\begin{description}
\setlength{\itemsep}{3pt}
\item[i.] Solve reduced eigenvalue problem $A' q=\lambda B' q$, with $A' = X'^T A X'$ and $B'=X'^TX'$
\item[ii.] Get new estimate for subspace $X$: $X_{i+1} = X' Q$   
\end{description}
\item[3.] Check the eigenvector error $r=max\  ||Ax_k-\lambda _k x_k||,\ 1 \leq k \leq m_0,\ \lambda_k \in {\cal I}$. If $r$ is above a given tolerance, GOTO 1.
\end{description}
}
\end{framed}
\end{minipage}
\medskip

\medskip

\noindent \begin{minipage}{\linewidth}
\begin{framed}
{\footnotesize
\noindent \textbf{XFEAST Algorithm} 

\noindent \hrulefill

\noindent {\bf Start with:} Matrix $A \in \mathbb{R}^{n \times n}$ to be diagonalized, 
interval ${\cal I}=(\lambda _{min},\lambda _{max} )$ where 
fewer than $m_0$ eigenvalues are expected to be found, initial guess $X_0 \in \mathbb{R}^{n\times m_0}$, and maximum number
of subspaces to store $s$. 
\begin{description}
\setlength{\itemsep}{3pt}
\item[0.] Repeatedly apply the filter procedure and append the resulting subspaces in order to 
expand the subspace to the predetermined size:
\begin{center} $X'=[X_0\ \ X_1\ \ X_2\ \ ...\ \ X_{s-1}]$ \end{center}
\noindent with $X_i=(\rho(A))^i X_0$.
\item[1.] Orthogonalize the columns of $X'$
\item[2.] Perform Rayleigh-Ritz procedure to find new estimate for eigenvalues and eigenvectors: 
\begin{description}
\setlength{\itemsep}{3pt}
\item[i.] Solve reduced eigenvalue problem $A' q=\lambda B' q$, with $A' = X'^T A X'$ and $B'=X'^TX'$
\item[ii.] Select the desired $m_0$ eigenpairs and get new estimate for $X$: $X_{i+1} = X' Q$
\end{description}
\item[3.] Check the eigenvector error $r=max\  ||Ax_k-\lambda _k x_k||,\ 1 \leq k \leq m_0,\ \lambda_k \in {\cal I}$. 
If $r$ is below a given tolerance, STOP.
\item[4.] Update subspace: 
\begin{description}
\setlength{\itemsep}{3pt}
\item[i.] Apply filter to new subspace estimate: $X_{i+1} = \rho (A) X_{i+1}$
\item[ii.] Update subspace by removing the oldest subspace and appending the newest update: $X'=[X_{i-s+1}\ \ X_{i-s+2}\ \ ...\ \ X_{i+1} ]$    
\end{description}
\item[5.] GOTO 1.

\end{description}
}
\end{framed}
\end{minipage}
\medskip

\noindent \begin{minipage}{\linewidth}
\begin{framed}
{\footnotesize
\noindent \textbf{RFEAST Algorithm} 

\noindent \hrulefill

\noindent {\bf Start with:} Matrix $A \in \mathbb{R}^{n \times n}$ to be diagonalized,
 interval ${\cal I}=(\lambda _{min},\lambda _{max} )$ where 
fewer than $m_0$ eigenvalues are expected to be found,
initial guess  $X_0 \in \mathbb{R}^{n\times m_0}$, maximum number of Rayleigh-Ritz iterations $s$ .
\begin{description}
\setlength{\itemsep}{3pt}
\item[1.] Filter the subspace $X_i$ to remove eigenvectors whose eigenvalues do not lie in the interval 
${\cal I}$: 
$X'=\rho (A)X_i$
\item[2.] Perform Rayleigh-Ritz procedure to find new estimate for eigenvalues and eigenvectors: 

\textbf{For } $j=1$ to $s$
\begin{description}
\setlength{\itemsep}{3pt}
\item[i.] Solve reduced eigenvalue problem $A' q=\lambda B' q$, with $A' = X'^T A X'$ and $B'=X'^TX'$
\item[ii.] Select the desired $m_0$ eigenpairs and get new estimate for $X$: $X_{i+1} = X' Q$   
\item[iii.] Compute residual vectors and expand subspace: \\$R=AX_{i+1}-X_{i+1}\Lambda \ \longrightarrow \ \ X'=[X' \ \ R]$
\end{description}
\textbf{end for}

\item[3.] Check the eigenvector error $r=max\  ||Ax_k-\lambda _k x_k||,\ 1 \leq k \leq m,\ \lambda_k \in {\cal I}$. If $r$ is above a given tolerance, GOTO 1.
\end{description}
}
\end{framed}
\end{minipage}

\end{document}